\documentclass[12pt,leqno]{article}
\usepackage{amsmath,amssymb,amsthm,amscd,dsfont}
\numberwithin{equation}{section} \allowdisplaybreaks
\begin{document}
\newtheorem{theorem}{Theorem}[section]
\newtheorem{defin}{Definition}[section]
\newtheorem{prop}{Proposition}[section]
\newtheorem{corol}{Corollary}[section]
\newtheorem{lemma}{Lemma}[section]
\newtheorem{rem}{Remark}[section]
\newtheorem{example}{Example}[section]
%\label{} %\ref{}
\title{From generalized K\"ahler to generalized Sasakian
structures}
\author{{\small by}\vspace{2mm}\\Izu Vaisman}
\date{}
\maketitle
{\def\thefootnote{*}\footnotetext[1]%
{{\it 2000 Mathematics Subject Classification: 53C15}.
Generalized Sasakian Structure.}}
\begin{center} \begin{minipage}{12cm}
A{\footnotesize BSTRACT. This is an expository paper, which provides
a first introduction to geometric structures on $TM\oplus T^*M$. The
paper contains definitions and characteristic properties of
generalized complex, generalized K\"ahler, generalized (normal,
almost) contact and generalized Sasakian structures. A few of these
properties are new.}
\end{minipage}
\end{center} \vspace*{5mm}
%\noindent
%begin{center} %\section %\end{center}
%\begin{center}
\section{Introduction} This is an expository paper.
Its aim is to introduce the reader into the new subject of
generalized structures. The non-previously published results are
Proposition \ref{bignormal} and Theorem \ref{thSgen}, which give
new characterizations of generalized, normal, almost contact and
generalized, Sasakian structures, and some remarks about non
degenerate, generalized, almost contact structures.

The word ``generalized" has the following precise meaning. If $M$
is an $m$-dimensional, differentiable manifold, a ``classical
structure" on $M$ is a reduction of the structure group of the
tangent bundle $TM$ from the general linear group
$Gl(m,\mathds{R})$ to a certain subgroup $G$. The
``generalization" consists in replacing $TM$ by the {\it big
tangent bundle} $T^{big}M=TM\oplus T^*(M)$. The bundle $T^{big}M$
has a natural, neutral metric (non degenerate, signature zero) $g$
defined by
\begin{equation}\label{g} g((X,\alpha),(Y,\beta))=
\frac{1}{2}(\alpha(Y)+\beta(X))\hspace{5mm}(X\in\chi(M),\alpha\in
\Omega^1(M)).
\end{equation} Hence, the natural structure group of $T^{big}M$ is
the group $O(m,m)$ that preserves the canonical neutral metric on
$\mathds{R}^{2m}$ and the ``generalized structures" will be
reductions of the structure group $O(m,m)$.

Furthermore, classical integrability conditions are expressed in
terms of the Lie bracket of vector fields on $M$. A generalized
bracket is the {\it Courant bracket}
\cite{C} given by the formula
\begin{equation}\label{crosetC}
[(X,\alpha),(Y,\beta)]=([X,Y],L_X\beta-L_Y\alpha
+\frac{1}{2}d(\alpha(Y)-\beta(X)).\end{equation} Accordingly,
generalized integrability conditions will be expressed in terms of
Courant brackets.

It is important that there are more natural bundle automorphisms
of $T^{big}M$ than those induced by diffeomorphisms of $M$. The
additional automorphisms are the $B$-{\it field transformations}
$(X,\alpha)\mapsto(X,\alpha+i(X)B)$ $(B\in\Omega^2(M),dB=0)$,
which also preserve the metric (\ref{g}) and the Courant bracket
(\ref{crosetC}).

The generalization described above is natural and appealing from
the differential geometric point of view. On the other hand, the
Courant bracket appeared in the study of constrained dynamical
systems and leads to important extensions of Hamiltonian mechanics
of a large interest in physics and control theory (see
\cite{VwH} and its references). Furthermore, the most
interesting generalized structures that were studied until now,
the generalized complex and the generalized K\"ahler structures,
appear in the study of supersymmetry in string theory (see
\cite{LRUZ} and its references). The study of generalized complex
structures as objects of differential geometry was initiated by
Hitchin \cite{Hitchin} and his student Marco Gualtieri, whose
thesis \cite{Gualt} turned into the standard reference on the
subject. Numerous authors, who also extended the scope of the
theory, have followed.
\section{Generalized complex and K\"ahler structures}
A generalized, almost complex structure is a reduction of the
structure group of $T^{big}M$ from $O(m,m)$ to $O(m,m)\cap
Gl(m,\mathds{C})$, therefore, it is defined by an endomorphism
$\mathcal{J}$ of $T^{big}M$ that has the properties
\begin{equation}\label{skewsym}
g(\mathcal{X},\mathcal{J}\mathcal{Y}) +
g(\mathcal{J}\mathcal{X},\mathcal{Y})=0,\;\;
\mathcal{J}^2=-Id,\end{equation} where calligraphic letters
denote pairs $\mathcal{X}=(X,\alpha)$, $\mathcal{Y}=(Y,\beta)$.
Equivalently, one may replace $\mathcal{J}$ by its
$i$-eigenbundle, which is a maximal (i.e.,
$rank_{\mathds{C}}L=m$), $g$-isotropic, complex subbundle $L$ of
the complexification $T^{big}M\otimes\mathds{C}$ such that
$L\cap\bar L=0$ (the bar denotes complex conjugation).

Furthermore, if the {\it Courant-Nijenhuis torsion} vanishes, i.e.,
\begin{equation}\label{tensNij}
\mathcal{N}_\mathcal{J}(\mathcal{X},\mathcal{Y}) =
[\mathcal{J}\mathcal{X},\mathcal{J}\mathcal{Y}] -
\mathcal{J}[\mathcal{X},\mathcal{J}\mathcal{Y}] -
\mathcal{J}[\mathcal{J}\mathcal{X},\mathcal{Y}]
+\mathcal{J}^2[\mathcal{X},\mathcal{Y}]=0,
\end{equation} $\forall \mathcal{X},\mathcal{Y}
\in\Gamma T^{big}M$ (the brackets are Courant brackets and $\Gamma$ denotes
the space of cross sections), the structure $\mathcal{J}$ is an
{\it integrable} or {\it generalized complex} structure and
$(M,\mathcal{J})$ is a generalized complex manifold. It follows
easily that integrability is equivalent with the fact that $\Gamma
L$ is closed under Courant brackets. In the terminology of
\cite{C}, adopted by all the authors in the field, $L$ is a
complex {\it Dirac structure}.

Similar definitions may be given for generalized, paracomplex
structures ($\mathcal{J}^2=Id$) and the study of the latter is
similar to that of the complex case \cite{Vgc}.

Generalized complex structures may be represented by classical
tensor fields as follows: \begin{equation}\label{matriceaJ}
\mathcal{J}\left(
\begin{array}{c}X\vspace{2mm}\\ \alpha \end{array}
\right) = \left(\begin{array}{cc} A&\sharp_\pi\vspace{2mm}\\
\flat_\sigma&-^t\hspace{-1pt}A\end{array}\right)
\left( \begin{array}{c}X\vspace{2mm}\\
\alpha \end{array}\right). \end{equation} Here, $A\in End(TM),
\pi\in\Gamma\wedge^2TM,\sigma\in\Omega^2(M)$, the musical
morphisms are defined like in Riemannian geometry and the index
$t$ denotes transposition. The first relation (\ref{skewsym})
explains the presence of $-^t\hspace{-1pt}A$ and the skew symmetry
of $\pi,\sigma$. Finally, the condition $\mathcal{J}^2=- Id$ is
equivalent with
\begin{equation}\label{2-reddezv} A^2=- Id -
\sharp_\pi\circ\flat_\sigma,\;\pi(\alpha\circ A,\beta)=\pi(\alpha,
\beta\circ A),\;\sigma(AX,Y)=\sigma(X,AY).\end{equation} The second,
respectively third, condition (\ref{2-reddezv}), are the {\it
compatibility} of $\pi$, respectively $\sigma$, with $A$.
\begin{rem}\label{obsdimpara} {\rm As a consequence of
(\ref{2-reddezv}), if a manifold $M$ has a generalized, almost
complex structure, the dimension of $M$ is even. Indeed, by
(\ref{2-reddezv})
$A|_{ker\,\flat_\sigma}:ker\,\flat_\sigma\rightarrow
ker\,\flat_\sigma$ and $A^2|_{ker\,\flat_\sigma}=- Id$. Thus,
$dim(ker\,\flat_\sigma)$ is even and, since
$dim(im\,\flat_\sigma)$ is even too, $dim M$ is even.}\end{rem}
Lengthy calculations show that the integrability of $\mathcal{J}$
given by (\ref{matriceaJ}) is equivalent with the following four
conditions \cite{{Cr},{Vgc}}:

i) the bivector field $\pi$ defines a Poisson structure on $M$
(i.e., the bracket of functions on $M$ given by
$\{f_1,f_2\}=\pi(df_1,df_2)$ is a Lie algebra bracket);

ii) the {\it Schouten concomitant}
\begin{equation}\label{SchoutenR}R(\pi,A)(\alpha,X) =
\sharp_\pi[L_X(\alpha\circ A)-L_{AX}\alpha]- (L_{\sharp_\pi\alpha}A)(X)
\end{equation} vanishes;

iii) the Nijenhuis tensor of $A$ (defined by (\ref{tensNij}) with
Lie brackets) satisfies the condition
\begin{equation}\label{Nijptintegrab} \mathcal{N}_A(X,Y) =
\sharp_\pi[i(Y)i(X)d\sigma];\end{equation}

iv) the {\it associated form}
\begin{equation}\label{sigmaA}\sigma_A(X,Y)=\sigma(AX,Y)\end{equation}
satisfies the condition
\begin{equation}\label{difsigmaA}
d\sigma_A(X,Y,Z)=\sum_{Cycl(X,Y,Z)}d\sigma(AX,Y,Z).\end{equation}

There are plenty of examples of interesting, generalized, complex
manifolds and we mention: (i) classical complex manifolds $(M,J)$,
with $A=J,\pi=0,\sigma=0$; (ii) classical symplectic manifolds
$(M,\omega)$, with $A=0,\pi=-\omega^{-1},\sigma=\omega$
$(d\omega=0)$; (iii) Hitchin pairs $(\varpi,A)$, where $\varpi$ is
a symplectic form, $A\in End(TM)$, $\varpi(AX,Y)=\varpi(X,AY)$,
and $\varpi_A(X,Y)=\varpi(AX,Y)$ is a closed $2$-form, with $\pi$
defined by $\flat_\varpi\circ\sharp_\pi=-Id$ and
$\sigma=\varpi\circ(Id+A^2)$ (this situation includes all the
generalized, complex structures with a non degenerate bivector
field $\pi$ \cite{Cr}); (iv) $\mathds{C}\mathds{P}^n$ and
manifolds obtained from $\mathds{C}\mathds{P}^2$ by blowing-up at
a finite number of points \cite{LT}, (v) the manifold
$M=3\mathds{C}\mathds{P}^2\#19\overline{\mathds{C}\mathds{P}^2}$,
which has neither a classical complex structure nor a classical
symplectic structure \cite{CG}. Notice also that any $B$-field
transformation sends a generalized complex structure to a new
generalized complex structure.

Now, we shall analyze the meaning of a {\it generalized Riemannian
structure} on $M$. Such a structure is a reduction of the
structure group of $T^{big}M$ from $O(m,m)$ to $O(m)\times O(m)$,
therefore, it consists of a $g$-orthogonal decomposition
$T^{big}M=V_+\oplus V_-$ where the terms are $m$-dimensional
subbundles and $g_+=g|_{V_+},g_-=g|_{V_-}$ are positive and
negative definite, respectively. This produces a positive definite
metric $G$ on $T^{big}M$ such that
\begin{equation}\label{G} \frac{1}{2}G=g_+-g_-.\end{equation}
The factor $1/2$ was introduced in order to ensure that $V_\pm$
(which are $G$-orthogonal) are the $\pm1$-eigenspaces of the
musical isomorphism
$$\sharp_G:T^*M\oplus TM\approx TM\oplus T^*M\rightarrow TM\oplus T^*M$$
defined by
\begin{equation}\label{defsharpG}
g(\sharp_G(X,\alpha),(Y,\beta))=\frac{1}{2}G((X,\alpha),(Y,\beta)).\end{equation}
(The isomorphism $\approx$ switches the order of the terms in a
pair.) Thus, a generalized, Riemannian structure may be seen as a
Euclidean (positive definite) metric $G$ on the bundle $T^{big}M$
such that $\sharp_G$ satisfies the conditions
\begin{equation}\label{Riemannbig}
\sharp_G^2=Id,\;\;g(\sharp_G(X,\alpha),\sharp_G(Y,\beta)) =
g((X,\alpha),(Y,\beta)).
\end{equation}
\begin{prop}\label{Gclasic} {\rm\cite{Gualt}} There exists a bijective
correspondence between generalized, Riemannian metrics $G$ and
pairs $(\gamma,\psi)$, where $\gamma$ is a classical Riemannian
metric and $\psi$ is a $2$-form on $M$.\end{prop}
\begin{proof} We may represent $G$ by
\begin{equation}\label{matriceaG} \sharp_G\left(\begin{array}{c}
X\vspace{2mm}\\ \alpha\end{array}\right)=
\left(\begin{array}{cc}\varphi&\sharp_\gamma\vspace{2mm}\\
\flat_\beta&^t\hspace{-1pt}\varphi\end{array}\right)
\left(\begin{array}{c} X\vspace{2mm}\\ \alpha\end{array}\right),
\end{equation} where $\varphi\in End(TM)$ and $\beta,\gamma$ are classical
Riemannian metrics on $M$ (use (\ref{Riemannbig}) and the fact
that $G$ is positive definite). Furthermore, the first condition
(\ref{Riemannbig}) is equivalent to
\begin{equation}\label{GcuRbig1} \begin{array}{c}
\varphi^2=Id-\sharp_\gamma\circ\flat_\beta,\;\gamma(\varphi X,Y)+
\gamma(X,\varphi Y)=0,\vspace{2mm}\\ \beta(\varphi
X,Y)+\beta(X,\varphi Y)=0.\end{array}\end{equation} Since $\gamma$
is non degenerate, the first condition (\ref{GcuRbig1}) yields
$\flat_\beta=\flat_{\gamma}\circ(Id-\varphi^2)$, hence, $G$
bijectively corresponds to the pair $(\gamma,\varphi)$. But,
$\varphi$ may be replaced by the $2$-form $\psi$ given by $
\flat_\psi=-\flat_\gamma\circ\varphi$, hence
$G
\leftrightarrow(\gamma,\psi)$ as claimed.\end{proof}

If $\psi=0$, we have $\varphi=0$, $\beta=\gamma$ and the
generalized metric reduces to a classical Riemannian metric. The
following result \cite{Gualt} is also important.
\begin{prop}\label{Vcupsi} The eigenbundles of the generalized,
Riemannian structure defined by the pair $(\gamma,\psi)$ are given
by the formula
\begin{equation}\label{exprVpm}
V_\pm=\{(X,\flat_{\psi\pm\gamma}X)\,/\,X\in TM\}.\end{equation}
\end{prop}
\begin{proof} The projectors associated with
the decomposition $T^{big}M=V_+\oplus V_-$ are given by
$pr_{\pm}=\frac{1}{2}(Id\pm\sharp_G)$ and, if we apply them to
pairs $(X,0)$, we see that the mappings
$\tau_{\pm}=pr_{TM}|_{V_\pm}$ are isomorphisms. Using
(\ref{matriceaG}) and the definition of $\psi$ we get
\begin{equation}\label{formuletau} \begin{array}{l}\tau_+^{-1}(X)=
(X,\flat_\gamma(X-\varphi
X))=(X,\flat_{\psi+\gamma}X),\vspace{2mm}\\
\tau_-^{-1}(X)= (X,-\flat_\gamma(X+\varphi
X))=(X,\flat_{\psi-\gamma}X),\end{array}\end{equation} whence
(\ref{exprVpm}).\end{proof}
\begin{rem}\label{remGpm} {\rm The isomorphisms $\tau_{\pm}$ may be
used to transfer the metrics $G|_{V_\pm}$ to metrics $G_\pm$ of
the tangent bundle $TM$ and the result is
$$
G_{\pm}(X,Y)=G(\tau_{\pm}^{-1}(X),\tau_{\pm}^{-1}(Y))= \pm
2g(\tau^{-1}_{\pm}(X),\tau^{-1}_{\pm}(Y))=2\gamma(X,Y).
$$}\end{rem}

A generalized almost complex structure $ \mathcal{J}$ is
compatible with a generalized metric $G$ if the structure group of
$T^{big}M$ is further reduced to $U(m/2)\times U(m/2)$.
Accordingly, this compatibility condition is
\begin{equation}\label{compatGJ} G(\mathcal{J}(X,\alpha),(Y,\beta))
+
G((X,\alpha),\mathcal{J}(Y,\beta))=0\,
\stackrel{(\ref{skewsym}),(\ref{defsharpG})}{\Leftrightarrow}\,
\sharp_G\circ\mathcal{J}=\mathcal{J}\circ\sharp_G
\end{equation} and, if (\ref{compatGJ}) holds, $(G,\mathcal{J})$
is a generalized, {\it almost Hermitian structure}; if $
\mathcal{J}$ is integrable, ``almost" is dropped.

The $(G,\mathcal{J})$-compatibility implies that
$(G,\mathcal{J}^c=\sharp_G\circ\mathcal{J})$ is again a
generalized, almost Hermitian structure, {\it complementary} to
$\mathcal{J}$, such that
$\mathcal{J}\circ\mathcal{J}^c=\mathcal{J}^c\circ\mathcal{J}$.
Conversely, if $(\mathcal{J},\mathcal{J}^c)$ is a commuting pair
of generalized, almost complex structures such that $G$ defined by
$\sharp_G=-\mathcal{J}\circ\mathcal{J}^c$ is positive definite,
then, $G$ is a generalized, Riemannian metric, which is compatible
with $
\mathcal{J}$ and $\mathcal{J}^c$. Many authors prefer to define a
generalized, almost Hermitian structure by the pair
$(\mathcal{J},\mathcal{J}^c)$.
\begin{prop}\label{prGaltK} {\rm\cite{Gualt}} A generalized, almost
Hermitian structure $(G,\mathcal{J})$ is equivalent with a
quadruple $(\gamma,\psi,J_+,J_-)$, where $(\gamma,J_{\pm})$ are
classical, almost Hermitian structures and
$\psi\in\Omega^2(M)$.\end{prop}
\begin{proof} The pair $(\gamma,\psi)$ is the one which
is equivalent with $G$. The $(G,\mathcal{J})$-compatibility
implies that $\mathcal{J}$ preserves the subbundles $V_\pm$,
hence, it induces endomorphisms $
\mathcal{J}_\pm\in End(V_\pm)$ such that $
\mathcal{J}_\pm^2=-Id.$ The latter can be transferred to the
almost complex structures
\begin{equation}\label{strJ}
J_{\pm}=\tau_{\pm}\circ\mathcal{J}_{\pm}\circ\tau_{\pm}^{-1}
\stackrel{(\ref{matriceaJ})}{=}A+\sharp_\pi\circ\flat_{\psi\pm\gamma}
\in
End\,TM,\end{equation} which are compatible with $\gamma$ because
$\gamma$ has been obtained by the similar transfer of
$G|_{V_\pm}$. Conversely, $(\gamma,\psi)$ define the subbundles
$V_\pm$ by (\ref{exprVpm}). The structures $J_\pm$ are transferred
by $\tau_\pm$ to structures $ \mathcal{J}_\pm$ on $V_\pm$ and $
\mathcal{J}=\mathcal{J}_++\mathcal{J}_-$ together with $G$ defined
by $(\gamma,\psi)$ is the corresponding generalized, almost
Hermitian structure.
\end{proof}
\begin{rem}\label{reconstuctieJ} {\rm The following relations
between $ \mathcal{J}$ and $J_\pm$ are also interesting. Formulas
(\ref{strJ}) yield
\begin{equation}\label{AdinJ} \begin{array}{c}
\sharp_\pi=\frac{1}{2}(J_+-J_-)\circ\sharp_\gamma,\vspace{2mm}\\
A=\frac{1}{2}[J_+\circ(Id-\sharp_\gamma\flat_\psi)+
J_-\circ(Id+\sharp_\gamma\flat_\psi)].\end{array}
\end{equation} The remaining entry of the matrix
representation of $ \mathcal{J}$ will be obtained from the matrix
expression of (\ref{compatGJ}), which yields among others
\begin{equation}\label{determinsigma} \varphi\circ
A+\sharp_\gamma\circ\flat_\sigma= A\circ
\varphi+\sharp_\pi\circ\flat_\beta\,\Leftrightarrow\,
\flat_\sigma=
\flat_\gamma\circ(A\circ Q-Q\circ A+\sharp_\pi\circ\flat_\beta).
\end{equation}}\end{rem}
\begin{example}\label{JJc} {\rm Assume
that the structure $(G,\mathcal{J})$ has the corresponding
quadruple $(\gamma,\psi,J_+,J_-)$. By its definition, the
complementary structure $\mathcal{J}^c$ satisfies the conditions
$\mathcal{J}^c_{\pm}=\pm\mathcal{J}_{\pm}$ and formula
(\ref{strJ}) shows that $(G,\mathcal{J}^c)$ has the corresponding
quadruple $(\gamma,\psi,J_+,-J_-)$. In the classical case,
$\psi=0$ and $J_-=-J_+$.}\end{example}

To continue our path towards generalized, K\"ahler manifolds, we
notice that, if the metric $G$ reduces to a classical Riemannian
metric $\gamma$ and the structure $\mathcal{J}$ reduces to a
classical structure $J$, the complementary structure
$\mathcal{J}^c$ is given by the matrix $$\mathcal{J}^c=
\left(
\begin{array}{cc} 0&\sharp_\pi\vspace{2mm}\\ \flat_\sigma&0
\end{array} \right)$$ where
$\sigma(X,Y)=\omega(X,Y)=\gamma(AX,Y),\;\pi=\sharp_\gamma\sigma$
($\omega$ is the {\it K\"ahler form} of $(\gamma,J)$). Thus, a
classical K\"ahler structure is characterized by the integrability
of the two structures $\mathcal{J}$ and $
\mathcal{J}^c$. Accordingly, one states \cite{Gualt}
\begin{defin}\label{defgenK} {\rm A {\it generalized K\"ahler
structure} is a generalized, almost Hermitian structure
$(G,\mathcal{J})$ such that the two structures
$\mathcal{J},\mathcal{J}^c$ are integrable. (We may also define a
{\it generalized, almost K\"ahler structure} $(G,\mathcal{J})$ by
requiring only the complementary structure $\mathcal{J}^c$ to be
integrable.)}\end{defin}

We shall prove the following theorem, which characterizes the
generalized K\"ahler structures.
\begin{theorem}\label{G4} {\rm\cite{VCRF}} The generalized almost
Hermitian structure $(G,\mathcal{J})$ with the associated
structures $(\gamma,\psi,J_\pm)$ is a generalized K\"ahler
structure iff $(\gamma,J_\pm)$ are classical Hermitian structures
and
\begin{equation}\label{nouG4} (\nabla_XJ_\pm)(Y)=\pm\frac{1}{2}
\sharp_\gamma[i(X)i(J_\pm Y)d\psi + (i(Y)i(X)d\psi)\circ J_\pm],
\end{equation} where $\nabla$ is the Levi-Civita connection of the
metric $\gamma$.\end{theorem}

This theorem is the consequence of a sequence of lemmas as
follows.
\begin{lemma}\label{G1} {\rm\cite{Gualt}} The generalized, almost Hermitian
structure $(G,\mathcal{J})$ is a generalized K\"ahler structure
iff the $i$-eigenbundles $L_\pm$ of $\mathcal{J}_\pm$ are closed
under Courant brackets.\end{lemma}
\begin{proof} Consider also the $i$-eigenbundles $L^c_\pm$ of
$\mathcal{J}^c_\pm$ and notice that the relations $
\mathcal{J}^c=\sharp_G\circ\mathcal{J},\sharp_G=-\mathcal{J}\circ\mathcal{J}^c$
imply
\begin{equation}\label{releigen}
\begin{array}{c}L=L_+\oplus L_-,\, L_+=L\cap V_+,\,L_-=L\cap
V_-,\,L^c_+=L_+,\vspace{2mm}\\L^c_-=\bar{L}^c_+, L_+=L\cap
L^c,\,L_-=L\cap\bar{L}^c.\end{array}\end{equation} Now, if the
structure $(G,\mathcal{J})$ is generalized K\"ahler, the last two
equalities (\ref{releigen}) obviously imply that $L_\pm$ are
closed under Courant brackets. To get the converse result, it is
enough to prove that $$ \mathcal{X}_+\in L_+,\mathcal{Y}_-\in L_-
\,\Rightarrow\,[ \mathcal{X}_+,\mathcal{Y}_-]\in L.$$
This follows from the following property of the Courant bracket
\cite{C} \begin{equation}\label{axv} \begin{array}{c}(pr_{TM}\mathcal{Z})(g(
\mathcal{X},\mathcal{Y}))=g([ \mathcal{Z},\mathcal{X}],\mathcal{Y})
+g(\mathcal{X},[ \mathcal{Z},\mathcal{Y}])\vspace{2mm}\\
+\frac{1}{2}(pr_{TM}
\mathcal{X})(g( \mathcal{Z},\mathcal{Y})) +\frac{1}{2}(pr_{TM}
\mathcal{Y})(g( \mathcal{Z},\mathcal{X})).
\end{array}\end{equation} If this equality is applied for
$(\mathcal{X},\mathcal{Y},\mathcal{Z})\mapsto
(\mathcal{Z}_+,\mathcal{Y}_-,\mathcal{X}_+)$ and
$(\mathcal{X},\mathcal{Y},\mathcal{Z})\mapsto
(\mathcal{X}_+,\mathcal{Z}_-,\mathcal{Y}_-)$, then, using the
$g$-orthogonality relations given by the isotropy of $L$, we get
$[\mathcal{X}_+,\mathcal{Y}_-]\perp_g\mathcal{Z}_\pm$, whence $[
\mathcal{X}_+,\mathcal{Y}_-]=0$, which is more than needed for the conclusion.
\end{proof}
\begin{lemma}\label{G2} The generalized, almost Hermitian structure
$(G,\mathcal{J})$ with the associated structures
$(\gamma,\psi,J_\pm)$ is a generalized K\"ahler structure iff
$(\gamma,J_\pm)$ are classical Hermitian structures and
\begin{equation}\label{eqG2} i(X\wedge Y)d\psi= \pm (i(X)L_Y\gamma
- L_Xi(Y)\gamma)\;\;(X,Y\in S_\pm),
\end{equation} where $S_\pm\subseteq TM$ are the
$i$-eigenbundles of  $J_\pm$. \end{lemma}
\begin{proof} Since $\mathcal{J}_\pm$ are the $\tau_\pm$-transfers
of $J_\pm$, we have
\begin{equation}\label{eqluiS} L_\pm=\{(X,\flat_{\psi+\gamma} X)\,/\,X\in
S_\pm\}= \{(X,\flat_{\psi-i\omega_\pm} X)\,/\,X\in
S_\pm\},\end{equation} where $\omega_\pm$ are the K\"ahler forms
of $(\gamma,J_\pm)$. Furthermore, we can get the following
expression of the required brackets
\begin{equation}\label{auxgeneral} [(X,\flat_{\psi\pm\gamma}X),
(Y,\flat_{\psi\pm\gamma}Y)] = ([X,Y],\flat_{\psi\pm\gamma}[X,Y]
\end{equation} $$+i(Y)i(X)d\psi
\pm (L_Xi(Y)\gamma-i(X)L_Y\gamma))\hspace{2mm}(X,Y\in\chi^1(M)).$$
This follows by evaluating the $1$-form component of the bracket
on a vector field $Z$. Formula (\ref{auxgeneral}) shows the
equivalence between the generalized K\"ahler conditions stated in
Lemmas
\ref{G1}, \ref{G2}.\end{proof}

Furthermore, take two unitary connections $\nabla^\pm$ on $TM$,
i.e., such that $\nabla^\pm\gamma=0,\; \nabla^\pm J_\pm=0$, and
consider the {\it difference tensors}
$\Theta^\pm(X,Y)=\nabla^\pm_XY-\nabla_XY,$ where $\nabla$ is the
Levi-Civita connection of the metric $\gamma$. From
$\nabla\gamma=0$, we get
\begin{equation}\label{condTheta} \gamma(\Theta^\pm(X,Y),Z)+
\gamma(Y,\Theta^\pm(X,Z))=0.\end{equation} On the other hand, the condition
$\nabla^\pm J_\pm=0$ is equivalent with
\begin{equation}\label{condTheta2}
\Theta^\pm(X,J_\pm Y)-J_\pm\Theta^\pm(X,Y)=-(\nabla_XJ_\pm)(Y).\end{equation}
\begin{lemma}\label{G3} Let $(G,\mathcal{J})$ be a
generalized Hermitian structure with the associated structures
$(\gamma,\psi,J_\pm)$ and let $\nabla^\pm$ be unitary connections.
Then, $(G,\mathcal{J})$ is a generalized K\"ahler structure iff
$(\gamma,J_\pm)$ are classical Hermitian structures and the
equalities
\begin{equation}\label{eqG3}
\gamma(\Theta^\pm(Z,Y),X)=\mp\frac{1}{2}
d\psi(X,Y,Z)\end{equation} hold for any $Z\in\chi^1(M)$ and any
$X,Y\in S_\pm$.\end{lemma}
\begin{proof}
By a simple computation that uses $\nabla^\pm\gamma=0$ we get
\begin{equation}\label{Liecutors} \begin{array}{r}(L_X\gamma)(Y,Z)=
\gamma(\nabla^\pm_YX,Z)
+\gamma(Y,\nabla^\pm_ZX)\vspace{2mm}\\ +\gamma(T^\pm(X,Y),Z)
+\gamma(Y,T^\pm(X,Z)),\end{array}\end{equation} where $T^\pm$ is
the torsion of $\nabla^\pm$. Then, if we evaluate (\ref{eqG2}) on
$Z\in\chi^1(M)$ and use (\ref{Liecutors}), we get the following
equivalent form of (\ref{eqG2}):
\begin{equation}\label{auxeqG2} \begin{array}{l}
d\psi(X,Y,Z)=\pm [\gamma(X,\nabla^\pm_ZY)-\gamma(Y,\nabla^\pm_ZX)
\vspace{2mm}\\ +\gamma(X,T^\pm(Y,Z)) +\gamma(Y,T^\pm(Z,X))
-\gamma(Z,T^\pm(X,Y)),\end{array}\end{equation} where the first
two terms of the right hand side vanish if $X,Y\in S_\pm$. If we
insert $$ T^\pm(X,Y)=
\Theta^\pm(X,Y)-\Theta^\pm(Y,X),$$
in (\ref{auxeqG2}), we get (\ref{eqG3}).
\end{proof}
\vspace*{2mm}
{\it Proof of Theorem \ref{G4}.} From (\ref{eqG3}), we get
conditions with general arguments $X,Y,Z\in\chi^1(M)$ by replacing
$X,Y\in S_\pm$ by $(Id-iJ_\pm)X,(Id-iJ_\pm)Y$. The resulting
conditions have a real and an imaginary part, which are equivalent
via the change $X\mapsto J_\pm X$. Thus, we remain with the
following characterization of the generalized K\"ahler structures
\begin{equation}\label{eqG5} \begin{array}{l}
\gamma(\Theta^\pm(Z,J_\pm Y), -X) +
\gamma(\Theta^\pm(Z,- Y), J_\pm X)\vspace{2mm}\\
=\mp\frac{1}{2}[d\psi (J_\pm X,J_\pm^2Y,Z) +d\psi (- X,J_\pm
Y,Z)],\vspace{2mm}\\
\end{array}\end{equation}
Now, if we use (\ref{condTheta2}) and the equality $$\nabla
J_\pm^2=0=J_\pm\circ\nabla J_\pm+\nabla J_\pm\circ J_\pm,$$ we get
the following system that is equivalent to (\ref{eqG5}):
\begin{equation}\label{eqG6}
\gamma(J_\pm X,(\nabla_ZJ_\pm)(J_\pm Y))=
\pm\frac{1}{2}[-d\psi(J_\pm X,Y,Z)+d\psi(-X,J_\pm
Y,Z)].\end{equation} This result is equivalent with the required
condition (\ref{nouG4}).\hspace{2.7cm}$\square$

Furthermore, one has
\begin{theorem}\label{G6} {\rm\cite{Gualt}}
The generalized almost Hermitian structure $(G,\mathcal{J})$ with
the associated structures $(\gamma,\psi,J_\pm)$ is a generalized
K\"ahler structure iff $(\gamma,J_\pm)$ are classical Hermitian
structures and
\begin{equation}\label{eqluiGaltieri}
d\omega_\pm(J_\pm X,J_\pm Y,J_\pm Z)=\pm d\psi(X,Y,Z),
\end{equation} where $\omega_\pm$ are the K\"ahler forms of
$(\gamma,J_\pm)$.\end{theorem}
\begin{proof} We already know that
$J_\pm$ are complex structures in the generalized K\"ahler case.
Accordingly, the following formula holds (e.g.,
\cite{KN})
$$ \gamma((\nabla_XJ_\pm)(Y),Z) =
\frac{1}{2}[d\omega_\pm(X,Y,Z)-d\omega_\pm(X,J_\pm Y,J\pm
Z)]$$ and (\ref{eqG6}) yields $$
\begin{array}{l}d\psi(X,Y,Z)-d\psi(J_\pm X,J_\pm
Y,Z)\vspace{2mm}\\=\pm[d\omega_\pm(J_\pm X,
Y,Z)+d\omega_\pm(X,J_\pm Y,Z)].\end{array}$$ Now, replace
$X\mapsto J_\pm X$, then, subtract the first cyclic permutation of
$(X,Y,Z)$ and add the second cyclic permutation. The result is
\begin{equation}\label{eqauxthG}\begin{array}{l}
d\psi(X,Y,Z)=\pm\frac{1}{2}[d\omega_\pm(J_\pm X,J_\pm Y,J_\pm
Z)\vspace{2mm}\\ + d\omega_\pm(J_\pm X,Y,Z)+d\omega_\pm(X,J_\pm
Y,Z)+d\omega_\pm(X,Y,J_\pm Z)].\end{array}
\end{equation}
Since for any Hermitian manifold $\omega_\pm$ is of the complex
type $(1,1)$ and $d\omega_\pm$ has no $(3,0)$ and $(0,3)$ type
components, if we use arguments in the eigenbundles of $J_\pm$, we
see that (\ref{eqauxthG}) coincides with
(\ref{eqluiGaltieri}).\end{proof}

The results above have analogous para-K\"ahler versions. On the
other hand the results were extended to generalized metric
$F$-structures $(F^3+F=0)$ \cite{VCRF}.

Concerning examples, $\mathds{C}\mathds{P}^n$ and
$\mathds{C}\mathds{P}^2$ blown-up at a finite number of points
were shown to be generalized K\"ahler in \cite{LT}. If the
$2$-form $\psi$ is closed, (\ref{nouG4}) reduces to $\nabla
J_\pm=0$, i.e., $(\gamma,J_\pm)$ are classical K\"ahler
structures. Therefore, any {\it bi-Hermitian manifold} $M$ is a
generalized K\"ahler manifold (add any closed $2$-form $\psi$ to
complete the structure). For instance, any hyper-K\"ahler manifold
has three bi-Hermitian structures. Bi-Hermitian manifolds were
studied and, in some cases, classified by several authors. The
reader will find more about the existence and non-existence of
generalized K\"ahler structures in \cite{AG} and its references.
\begin{center}
\section{Generalized almost contact and Sasakian
structures}
\end{center}
In the realm of classical structures, odd-dimensional
correspondents of complex and K\"ahler structures are obtained by
using complex and K\"ahler structures on the manifold
$M\times\mathds{R}$ as follows.

Let $J$ be a complex structure on $M^{2n+1}\times\mathds{R}$ such
that (i) $J$ is invariant by translation along $\mathds{R}$ and
(ii) $J(T\mathds{R})\subseteq TM$. Then $J$ is said to be
{\it$M$-adapted}. If we denote by $t$ the coordinate on
$\mathds{R}$, $J$ is an $M$-adapted structure iff
\begin{equation}\label{Jacontact} J=F+dt\otimes
Z-\xi\otimes\frac{\partial}{\partial t},\end{equation} where $F\in
End(TM), Z\in\chi(M),\xi\in\Omega^1(M).$ Accordingly, condition
$J^2=-Id$ becomes
\begin{equation}\label{almct} F^2=-Id+\xi\otimes Z,\;\;\xi\circ
F=0,\;FZ=0,\;\xi(Z)=1,\end{equation} and the triple $(F,Z,\xi)$ is
called an {\it almost contact structure} on $M$; it corresponds to
a reduction of the structure group of $TM$ to
$Gl(n,\mathds{C})\times\{1\}$. If the adapted structure $J$ is
integrable, the almost contact structure $(F,Z,\xi)$ is {\it
normal} and the normality condition is \cite{Bl}
\begin{equation}\label{normalact} \mathcal{N}_F+Z\otimes d\xi=0,
\end{equation} where $\mathcal{N}_F$ is the Nijenhuis tensor of
$F$.

A further reduction of the structure group of $TM$ to
$U(n)\times\{1\}$ is obtained by adding a Riemannian metric
$\gamma$ on $M$ such that the translation invariant, almost
complex structure $J$, is Hermitian with respect to the metric
\begin{equation}\label{gamma}
\Gamma=e^{t}(\gamma+dt^2)\end{equation} (the factor $e^t$ is
needed for the Sasakian structures, which will be defined later
on).

Then $J$ necessarily is $M$-adapted and we get an {\it almost
contact metric structure} $(F,Z,\xi,\gamma)$, where
\begin{equation}\label{acmetric}
\gamma(FX,FY)=\gamma(X,Y)-\xi(X)\xi(Y),
\end{equation}
which also implies $\xi=\flat_\gamma Z,g(Z,FX)=0,g(Z,Z)=1$
\cite{Bl}.

The almost contact metric structure $(F,Z,\xi,\gamma)$ has the
associated {\it fundamental $2$-form} $\Xi(X,Y)=g(FX,Y)$, while
the corresponding almost Hermitian structure $J$ has the K\"ahler
form $\omega$. A simple calculation gives
\begin{equation}\label{formaomega}
\omega=e^{t}(\Xi-\xi\wedge
dt),\;d\omega=e^{t}[d\Xi+(\Xi-d\xi)\wedge dt].\end{equation}

The most usual definition of a Sasakian structure requires it to
be a normal, contact, metric structure $(F,Z,\xi,\gamma)$ where
the use of the term {\it contact} instead of {\it almost contact}
means the requirement $\Xi=d\xi$. From (\ref{formaomega}), we see
that a Sasakian structure is characterized by the fact that the
corresponding structure $(\Gamma,J)$ is K\"ahler \cite{Bl}.

We will give generalized versions of the classical structures
recalled above.
\begin{defin}\label{bigadapt} {\rm A generalized, almost complex structure
$\mathcal{J}$ on $M\times\mathds{R}$ is said to be $M$-{\it
adapted} if it has the following three properties (a) $
\mathcal{J} $ is invariant by translation along $\mathds{R}$, (b)
$ \mathcal{J} (T\mathds{R}\oplus0)\subseteq 0\oplus T^*M$, (c)
$\mathcal{J}(0\oplus T^*\mathds{R})\subseteq
TM\oplus0$.}\end{defin}

The invariance of $ \mathcal{J}$ by translations means that the
Lie derivatives $L_{\partial/\partial t}$ of the classical tensor
fields of $
\mathcal{J}$ (defined by (\ref{matriceaJ})) vanish.
If conditions (b), (c) are also imposed, it follows that the
classical tensor fields of an $M$-adapted, generalized, almost
complex structure are of the form
\begin{equation}\label{genadapt} A= F,\,
\pi=P+Z\wedge\frac{\partial}{\partial t},\,\sigma=\theta+\xi\wedge
dt,\end{equation} where
$P\in\Gamma\wedge^2TM,\theta\in\Omega^2(M),Z\in\chi(M),\xi\in\Omega^1(M)$
(one may use local coordinates $x^i$ on $M$ to see this).
\begin{rem}\label{remtrinv} {\rm If $
\mathcal{J}$ only is translation invariant, the second and third formula
(\ref{genadapt}) hold but $A:TM\oplus\mathds{R}\rightarrow
TM\oplus\mathds{R}$ also includes a vector field
$V=pr_{TM}A(\partial/\partial t)$, a $1$-form
$\kappa(X)=pr_{\mathds{R}}(AX)$ and a function $s$ given by
$pr_{\mathds{R}}(A(\partial/\partial t))=s(\partial/\partial t)$.}
\end{rem} Furthermore, conditions (\ref{skewsym}) are equivalent to
\begin{equation}\label{condF}\begin{array}{l}
F\circ\sharp_P=\sharp_P\circ\hspace*{1pt}^t\hspace*{-1pt} F,\;
\flat_\theta\circ F=\hspace*{1pt}^t\hspace*{-1pt}
F\circ\flat_\theta,\;i(Z)\theta=0,\;i(\xi)P=0,\vspace{2mm}\\
F^2=-Id-\sharp_P\circ\flat_\theta+\xi\otimes Z,\;F(Z)=0,\;\xi\circ
F=0,\;\xi(Z)=1.
\end{array}\end{equation}

The triple $(F,P,\theta)$ defines an endomorphism $ \mathcal{F}$ of
$T^{big}M$ of matrix form \begin{equation}\label{matriceaF}
\mathcal{F}\left(
\begin{array}{c}X\vspace{2mm}\\ \alpha \end{array}
\right) = \left(\begin{array}{cc} F&\sharp_P\vspace{2mm}\\
\flat_\theta&-^t\hspace{-1pt}F\end{array}\right)
\left( \begin{array}{c}X\vspace{2mm}\\
\alpha \end{array}\right). \end{equation}
The pair $(Z,\xi)$ defines the endomorphism $ \mathcal{Z}$ of
$T^{big}M$ of matrix form
\begin{equation}\label{matriceaZ} \mathcal{Z}\left(
\begin{array}{c}X\vspace{2mm}\\ \alpha \end{array}
\right) = \left(\begin{array}{cc} Z\otimes\xi&0\vspace{2mm}\\
0&^t\hspace{-1pt}(Z\otimes\xi)\end{array}\right)
\left( \begin{array}{c}X\vspace{2mm}\\
\alpha \end{array}\right), \end{equation} where
$Z\otimes\xi:TM\rightarrow TM$ is the evaluation of $\xi$ and
$^t\hspace{-1pt}(Z\otimes\xi):T^*M\rightarrow T^*M$ is the
evaluation of $Z$. The conditions (\ref{condF}) are equivalent to
\begin{equation}\label{gac2} \flat_g\circ\mathcal{F}
+\hspace{+1pt}^t\hspace{-1pt}\mathcal{F}\circ\flat_g=0,\;\mathcal{F}^2=-Id+
\mathcal{Z},\; \mathcal{F}\circ
\mathcal{Z}=0,\;\|Z\oplus\xi\|_g=1.\end{equation} (The first
condition (\ref{gac2}) ensures that $P$ and $\theta$ are skew
symmetric, and $g$ is the neutral metric of $T^{big}M$.)
Accordingly, we define
\cite{PW},
\cite{VCRF}
\begin{defin}\label{genalmct} {\rm A {\it generalized almost contact
structure} on $M$ is a couple $(\mathcal{F}\in End(T^{big}M),
(Z,\xi)\in\Gamma T^{big}M)$ that satisfies (\ref{gac2}).
Equivalently the structure is a system of classical tensor fields
$(P,\theta,F,Z,\xi)$ that satisfies (\ref{condF}).}\end{defin}

We mention a few examples \cite{PW}. If $(F,Z,\xi)$ is a classical
almost contact structure, then, $(F,P=0,\theta=0,Z,\xi)$ is a
generalized, almost contact structure. If $\xi$
($\xi\wedge(d\xi)^n|_x\neq0$, $\forall x\in M$) is a contact form,
if $Z$ is the corresponding Reeb vector field (which means that
$\xi(Z)=1$ and $i(Z)d\xi=0$) and if $\theta=d\xi$, then,
$\phi(X)=\flat_{\theta}(X)-\xi(X)\xi$ is an isomorphism
$TM\rightarrow T^*M$, and we get a bivector field
$P(\alpha,\beta)=\theta(\phi^{-1}\alpha,\phi^{-1}\beta)$. Then,
$(F=0,P,\theta,Z,\xi$) is  a generalized, almost contact
structure. Thus, while a contact form has no canonically
associated, classical, almost contact structure, it defines a
canonical generalized, almost contact structure. A similar
situation holds for an almost cosymplectic structure
$(\xi\in\Omega^1(M),\theta\in\Omega^2(M))$ where
$\xi\wedge\theta^n$ nowhere vanishes.

A generalized, almost contact structure will be called {\it
normal} if the corresponding $M$-adapted, generalized, almost
complex structure on $M\times\mathds{R}$ is integrable. Thus, the
normality conditions are conditions i)-iv) of Section 1 applied to
the tensor fields (\ref{genadapt}). After some technical efforts,
it turns out that the normality conditions are:
\begin{equation}\label{normalitate} \begin{array}{l}
[P,P]=0,\:R_{(P,F)}=0,\vspace{2mm}\\ L_{Z}P=0,\:
L_{Z}\theta=0,\:L_{\sharp_P\alpha}\xi=0,
\vspace{2mm}\\
\mathcal{N}_F(X,Y)=\sharp_P(i(X\wedge Y)d\theta) - (d\xi(X,Y))Z,\vspace{2mm}\\
d\theta_F(X_1,X_2,X_3)=\sum_{Cycl(1,2,3)}d\theta(FX_1,X_2,X_3),\vspace{2mm}\\
L_{Z}\xi=0,\;L_{Z}F=0,\; (L_{FX}\xi)(Y) -
(L_{FY}\xi)(X)=0,\end{array}
\end{equation}  and, if at no point has $\sharp_P\circ\flat_\theta$
the eigenvalue $-1$, the last three conditions (\ref{normalitate})
follow from the other conditions \cite{Vstable}.

It is also possible to characterize normality by properties of the
couple $( \mathcal{F},(Z,\xi))$:
\begin{prop}\label{bignormal} The generalized, almost contact
structure $( \mathcal{F},(Z,\xi))$ is normal iff the following
conditions hold
\begin{equation}\label{eqbignormal} \begin{array}{l}
\mathcal{N}_{\mathcal{F}}((X,\alpha),(Y,\beta))=\mathcal{Z}([(X,\alpha),(Y,\beta)]),
\vspace{2mm}\\

[(Z,0),\mathcal{F}(X,\alpha)]=\mathcal{F}(L_ZX,L_Z\alpha),
\;\;[\mathcal{F}(X,\alpha),(0,\xi)]=\mathcal{F}(0,L_X\xi),
\end{array} \end{equation}
where the brackets are Courant brackets, in the first condition
$(X,\alpha),(Y,\beta)\in\Gamma T^{big}M$, and in the second and
third condition $(X,\alpha)\in im\,\mathcal{F}$.
\end{prop}
\begin{proof} Conditions (\ref{gac2}) imply $(X,\alpha)\in
im\,\mathcal{F}$ iff $\xi(X)=0,\alpha(Z)=0$, whence,
\begin{equation}\label{descmare}\begin{array}{c}
T^{big}(M\times\mathds{R})=im\mathcal{F}\oplus span\{(Z,0)\}\oplus
span\{(0,\xi)\}\vspace{2mm}\\ \oplus
span\{(\frac{\partial}{\partial t},0)\}\oplus span\{(0,dt)\}.
\end{array}\end{equation} On the other hand, from (\ref{genadapt})
we get $$ \mathcal{J}(X+f\frac{\partial}{\partial
t},\alpha+\varphi dt) $$
$$=(FX+\sharp_P\alpha+\alpha(Z)\frac{\partial}{\partial t}-\varphi
Z,\flat_\theta-\alpha\circ F+\xi(X)dt-f\xi)$$
$$=\mathcal{F}(X,\alpha)-\varphi(Z,0)-f(0,\xi) +
(\alpha(Z)\frac{\partial}{\partial t},\xi(X)dt),$$ whence,
$$\mathcal{J}(Z,0)=(0,dt),\;\mathcal{J}(0,\xi)=(\frac{\partial}{\partial t},0),\;
\mathcal{J}(X,\alpha)=\mathcal{F}(X,\alpha)\;\forall
(X,\alpha)\in im\,\mathcal{F}.$$ Now, we get the normality
conditions by asking $
\mathcal{N}_{\mathcal{J}}$ to vanish for all the possible
combinations of arguments in the various terms of decomposition
(\ref{descmare}). Moreover, since $
\mathcal{N}_{\mathcal{J}}$ is a tensor on $T^{big}M$, we do not
have to consider tensorial coefficients and just take arguments
$(X,\alpha)\in im\,\mathcal{F}$ and $(Z,0),(0,\xi),
(\frac{\partial}{\partial t},0),(0,dt)$.

The arguments $(\frac{\partial}{\partial t},0),(0,dt)$ produce the
important condition $L_Z\xi=0$, which is equivalent with the first
condition (\ref{eqbignormal}) for the arguments
$(X,\alpha)=(Z,0),(Y,\beta)=(0,\xi)$.

If this condition is used, the arguments $(X,\alpha)\in
im\,\mathcal{F},(0,dt)$ yield the second condition
(\ref{eqbignormal}), the arguments $(X,\alpha)\in
im\,\mathcal{F},(\frac{\partial}{\partial t},0)$ yield the third
condition and the arguments $(X,\alpha), (Y,\beta)\in
im\,\mathcal{F}$ yield the first condition (\ref{eqbignormal}) for
this situation.

Finally, if we consider the arguments $(X,\alpha)\in
im\,\mathcal{F},(Z,0)$, $(X,\alpha)\in im\,\mathcal{F},(0,\xi)$,
respectively, we get the first condition (\ref{eqbignormal}) for
these cases and the supplementary equality
$\xi([Z,FX+\sharp_P\alpha])=0$. The latter is satisfied since from
(\ref{condF}) and $L_Z\xi=0$ we have
$$i([Z,FX+\sharp_P\alpha])\xi=L_Zi(FX+\sharp_P\alpha)\xi
- i(FX+\sharp_P\alpha)L_Z\xi=0.$$

Other choices of the arguments do not lead to new
conditions.\end{proof}

We will say that a generalized, almost contact structure is
non-degenerate if the corresponding structure $\mathcal{J}$ of
Definition \ref{bigadapt} is non degenerate, i.e., the bivector
field $\pi$ given by (\ref{genadapt}) is non degenerate.
Equivalently, this means that $Z\wedge P^n\neq0$ at every point
$x\in M$, hence $TM=im\,\sharp_P\oplus span\{Z\}$. The
corresponding Hitchin pair (see the examples of generalized,
complex structures in Section 1 or \cite{Cr}) is $(\varpi,F)$ with
$F$ given by (\ref{genadapt}) and
$$\varpi=\omega+\xi\wedge dt$$ where $\omega\in\Omega^2(M)$ is the
unique $2$-form that satisfies the conditions
$$i(\sharp_P\lambda)\omega=-\lambda+\lambda(Z)\xi,i(Z)\omega=0.$$
The $F$-compatibility condition $\varpi(FX,Y)=\varpi(X,FY)$ is
equivalent with
$$\omega(FX,Y)=\omega(X,FY),\,\xi\circ F=0.$$ Thus, $\omega_F(X,Y)
=\omega(FX,Y)$ is a $2$-form, and we have $\xi\wedge\omega^n\neq0$
at every point of $M$. (Use (\ref{condF}) while checking all the
above.)

Then, the structure is normal iff
$\varpi,\varpi_F\in\Omega^2(M\times\mathds{R})$ are closed, which
is equivalent to
$$d\xi=0,d\omega=0,d\omega_F=0.$$
Hence, a non degenerate, generalized, almost contact structure is
equivalent with an almost cosymplectic structure
$\xi\in\Omega^1(M),\omega\in\Omega^2(M)$
($\xi\wedge\omega^n\neq0$) complemented by $F\in End(TM)$, which
is compatible with $\omega$ and such that $\xi\circ F=0$. If these
tensor fields are given, we get $Z$ by asking
$i(Z)\xi=1,i(Z)\omega=0$, we get $P$ from
$\flat_\varpi\circ\sharp_\pi=-Id$ and we get $\theta$ from the
conditions (\ref{condF}). In the case of a normal structure, the
almost cosymplectic structure is cosymplectic, i.e.,
$d\xi=0,d\omega=0$ and we also have $d\omega_F=0$.

If instead of normality we require the generalized almost complex
structure $\mathcal{J}'$ on $M\times\mathds{R}$ with the classical
tensor fields
$$ A= F,\,
\pi=e^t(P+Z\wedge\frac{\partial}{\partial
t}),\,\sigma=e^{-t}(\theta+\xi\wedge dt)$$ (obtained by a conformal
change of $\mathcal{J}$ in the sense of \cite{{Vstable},{VIasi}}) to
be integrable, then the generalized structure is equivalent with the
complemented, almost cosymplectic structure $(\xi,\omega=d\xi,F)$
where $\xi$ is a contact form on $M$. This observation shows that
the generalized, almost contact structures with the property that
$\mathcal{J}'$ is integrable (but, need not be non degenerate)
deserve to be called {\it generalized contact structures}. The
integrability of $\mathcal{J}'$ is equivalent with the fact that the
restriction of its $i$-eigenbundle to $t=0$ is a
$\mathcal{E}^1$-Dirac (Dirac-Jacobi, stable Dirac-Jacobi) structure
\cite{{IW},{Vstable}}. These structures are integrable, generalized,
almost contact in the sense of \cite{IW}; however, the latter is a
larger class of structures since the corresponding structure
$\mathcal{J}$ is only required to be translation invariant.

The corresponding integrability conditions can be derived from the
integrability conditions of a generalized, almost complex structure
given in Section 1 (it is convenient to use Proposition 3.1 of
\cite{VIasi} as an intermediary step) and the results are equivalent
to
\begin{equation}\label{contactgen} \begin{array}{l}
[P,P]=2Z\wedge P,\:R_{(P,F)}=0,\vspace{2mm}\\ L_{Z}P=0,\:
L_{Z}\theta=0,\:L_{\sharp_P\alpha}\xi=\flat_\theta\sharp_P\alpha,
\vspace{2mm}\\
\mathcal{N}_F(X,Y)=\sharp_P(i(X\wedge Y)d\theta) - (d\xi(X,Y)-\theta(X,Y))Z,\vspace{2mm}\\
d\theta_F(X_1,X_2,X_3)=\sum_{Cycl(1,2,3)}d\theta(FX_1,X_2,X_3),\vspace{2mm}\\
L_{Z}\xi=0,\;L_{Z}F=0,\; (L_{FX}\xi)(Y) -
(L_{FY}\xi)(X)=\theta_F(X,Y).\end{array}
\end{equation} In particular, the tensor fields $(P,Z)$ define a Jacobi structure
on $M$. By comparing (\ref{contactgen}) with (\ref{normalitate})
we see that a generalized contact structure in this sense is
normal iff $P=0,\theta=0$ and $(F,Z,\xi)$ is a classical, normal,
almost contact structure. On the other hand, if $F=0$, since $Z$
is not in the image of $\sharp_P$, the Nijenhuis tensor condition
in (\ref{contactgen}) implies $\theta=d\xi$ and the structure
reduces to that associated to a contact form and its Reeb vector
field.

A different terminology is proposed in \cite{PW} by the
introduction of two other notions. Namely, since a generalized,
almost contact structure $\mathcal{F}$ satisfies $
\mathcal{F}^3+\mathcal{F}=0$ it has the eigenvalues $\pm i,0$ and
corresponding eigenbundles $E_\pm, S\subseteq T^{big}M$ where
$S=span\{(Z,0),(0,\xi)\}$. Denote $L=E_+\oplus span\{(Z,0)\},
L^*=E_-\oplus span\{(0,\xi)\}$. In
\cite {PW}, $\mathcal{F}$ is a {\it generalized, contact structure} if $L$ is
closed under Courant brackets and a {\it strong, generalized,
contact structure} if both $L$ and $L^*$ are closed under Courant
brackets. The names were chosen because $L$ is bracket-closed in
the case of a contact form, while $L^*$ is not. On the other hand,
one has the ``strong contact" situation in the case of a
cosymplectic structure. It is easy to see that the $\pm
i$-eigenbundles $T_\pm$ of the corresponding, generalized, almost
complex structure $
\mathcal{J}$ of $M\times \mathds{R}$ are given by
$$T_\pm=E_\pm\oplus
span\{(Z,0)\mp i(0,dt),(0,\xi)\mp i(\frac{\partial}{\partial
t},0)\}.$$ With this formula, we can check that a generalized,
almost contact structure is normal iff it is a strong,
generalized, contact structure that satisfies the condition
$L_Z\xi=0$. All the examples of strong, generalized, contact
structures given in
\cite{PW} (cosymplectic manifolds, the $3$-dimensional Heisenberg
group) satisfy the condition $L_Z\xi=0$, hence, also are examples
of normal, generalized, almost contact structures.

The normal, generalized, almost contact manifolds
$(M,P,\theta,F,Z,\xi)$ have a nice geometric structure, which we
have described in \cite{Vstable}.
\begin{theorem}\label{strgnormalac} A generalized, normal, almost
contact structure on $M$ is equivalent with the following system
of geometric objects: 1) a vector field $Z$ whose trajectories
define a one-dimensional foliation $\mathfrak{Z}$, 2) a
complementary subbundle $\nu\mathfrak{Z}$ of $T\mathfrak{Z}$
$(T\mathfrak{Z}\oplus\nu\mathfrak{Z}=TM)$, 3) a transversal,
projectable, generalized, complex structure $
\mathfrak{J}$ of $\mathfrak{Z}$ with corresponding classical
tensor fields $F\in
End(\nu\mathfrak{Z}),P\in\Gamma\wedge^2\nu\mathfrak{Z},\theta\in
\Gamma\wedge^2(ann\,\mathfrak{Z})$,
such that the following properties hold: (i) $\nu\mathfrak{Z}$ is
invariant by the infinitesimal transformations $Z,\sharp_P\lambda$
$(\forall\lambda\in ann\,\mathfrak{Z})$, (ii) the Ehresmann
curvature of $\nu\mathfrak{Z}$ is invariant by $F$.\end{theorem}
\begin{proof} If we start with the tensor fields
$P,\theta,F,Z,\xi$, since $\xi(Z)=1$, $Z$ never vanishes and
defines a foliation $ \mathfrak{Z}$. A complementary bundle is
defined by $\nu\mathfrak{Z}=ann\,\xi$. By restricting
$(F,P,\theta)$ to $\nu\mathfrak{Z},\nu^*\mathfrak{Z}=ann\,Z$ we
get a generalized, almost complex structure $ \mathfrak{J}$ on
$\nu\mathfrak{Z}$ and by its being complex and projectable we
understand that it is projection-related with generalized, complex
structures on local transversal submanifolds of $ \mathfrak{Z}$.
This property of $
\mathfrak{J}$ and properties (i), (ii) follow from
the normality conditions (\ref{normalitate}). The details are
lengthy and we refer the interested reader to \cite{Vstable},
where the theorem is proven for more general structures ``of
codimension $h$" (i.e., with $h$ commuting vector fields $Z_h$).
We only recall that the Ehresmann curvature is defined by
$$R_{\nu\mathfrak{Z}}(X,Y)=-pr_{T\mathfrak{Z}}[pr_{\nu\mathfrak{Z}}X,
pr_{\nu\mathfrak{Z}}Y]$$ and its $F$-invariance means
$$R_{\nu\mathfrak{Z}}(FX,FY)=R_{\nu\mathfrak{Z}}(X,Y).$$

Conversely, if we start with $ Z,\nu\mathfrak{Z},\mathfrak{J}$ with
the required properties, we get a $1$-form $\xi$ by asking that
$\xi(Z)=1,\xi|_{\nu\mathfrak{Z}}=0$, and we have $
\mathfrak{Z}$-adapted local coordinates $(z,y^u)$ (i.e., $
\mathfrak{Z}$ is $y^u=0$ and $Z=\partial/\partial z$) such that
$\xi=dz+\xi_udy^u$ and $\nu\mathfrak{Z}$ has local bases
$Y_u=\partial/\partial y^u-\xi_u(\partial/\partial z)$. Then, the
tensor fields of $ \mathfrak{J}$ will be of the form $$
P=\frac{1}{2}P^{uv}(y^w)Y_u\wedge
Y_v,\;\theta=\frac{1}{2}\theta_{uv}(y^w)dy^u\wedge
dy^v,\;F(Y_u)=F_u^v(y^w)Y_v.$$ Again, a careful comparison between
properties (i), (ii) and the normality conditions
(\ref{normalitate}) shows that $(F,P,\theta,Z,\xi)$ is a normal,
generalized, almost contact structure on $M$
\cite{Vstable}.\end{proof}
\begin{example}\label{exnormal} {\rm The total space of a flat
principal circle bundle over a generalized, complex manifold is a
normal, generalized, almost contact manifold. Namely, with the
notation of Theorem \ref{strgnormalac}, we will take $Z$ to be the
fundamental, vertical vector field and $\xi$ to be the flat
connection form, then, $\nu\mathfrak{Z}$ will be given by $\xi=0$
and $P,F,\theta$ will be the lifts of the tensor fields of the
generalized, complex structure on the basis.}
\end{example}
\begin{rem}\label{transvingcont} {\rm Conditions
(\ref{contactgen}) similarly show that what we called a
generalized, contact manifold also has the foliation
$\mathfrak{Z}$, its transversal distribution $\nu\mathfrak{Z}$ and
the transversal, generalized, complex structure given by the
projections of $(F,P,\theta)$ but, we do not have a nice,
corresponding, characterization of the whole structure.}\end{rem}

Now, we shall bring a metric into the picture. With the classical
case in mind, we have to endow $M$ with a generalized Riemannian
metric $G$, equivalently, with a pair $(\gamma,\psi)$, where
$\gamma$ is a classical Riemannian metric and
$\psi\in\Omega^2(M)$, and associate to it a generalized,
Riemannian metric $\tilde G$ of $M\times\mathds{R}$ defined by a
pair
\begin{equation}\label{tildeG} \tilde G\Leftrightarrow
(\Gamma=e^{t}(\gamma+dt^2),\Psi=e^{t}(\psi+\kappa\wedge
dt)),\end{equation} where $\kappa\in\Omega^1(M)$. We skip a
discussion of {\it generalized, almost contact, metric
structures}, which seem to be less interesting, and directly go to
generalized, Sasakian structures
\cite{VCRF}.\begin{defin}\label{defgenS} {\rm A {\it generalized,
Sasakian manifold} is a generalized Riemannian manifold $(M,G)$
endowed with a translation invariant, generalized, almost complex
structure $\mathcal{J}$ of $M\times\mathds{R}$ such that, for some
$\kappa\in\Omega^1(M)$, $(M\times\mathds{R},\tilde G,\mathcal{J})$
is a generalized, K\"ahler manifold.}\end{defin}

Generalized, Sasakian manifolds exist. Indeed, if we define a {\it
bi-Sasakian structure} as a pair of different, classical, Sasakian
structures with the same metric, we get
\begin{prop}\label{existSgen} A bi-Sasakian structure on $M$, supplemented by
a $1$-form $\kappa\in\Omega^1(M)$, defines a generalized Sasakian
structure.\end{prop}
\begin{proof} Since a Sasakian structure on $(M,\gamma)$ is equivalent with a
K\"ahler structure on $(M\times\mathds{R},\Gamma)$, a given
bi-Sasakian structure $(F_\pm,Z_\pm,\xi_\pm,\gamma)$ is equivalent
with a bi-Hermitian structure on $(M\times\mathds{R},\Gamma)$ and,
if the latter is supplemented by a closed form
$\Psi=e^t(\psi+\kappa\wedge dt)$, a generalized K\"ahler
structure, i.e., a generalized Sasakian structure on $M$, will
arise. But, $d\Psi=0$ is equivalent with $d\psi=0,\psi+d\kappa=0$,
therefore, a choice of $\kappa$ will fix the generalized Sasakian
structure.\end{proof}

The generalized, Sasakian structures of Proposition \ref{existSgen}
will be said to be of bi-Sasakian type.
\begin{corol}\label{3S} If $(M,(F_\pm,Z_\pm,\xi_\pm,\gamma),\kappa)$
is a generalized Sasakian structure of bi-Sasakian type, then, one
of the following three situations occur: (i)
$F_-=-F_+,Z_-=-Z_+,\xi_-=-\xi_+$ (i.e., $M$ is a Sasakian manifold
with the two conjugated Sasakian structures), (ii) the structures
$(F_\pm,Z_\pm,\xi_\pm,\gamma)$ belong to a $3$-Sasakian structure,
(iii) the metric $\gamma$ must be of constant sectional curvature
$1$.\end{corol}
\begin{proof} These are known results for bi-Sasakian structures,
e.g., Lemmas 8.1.16, 8.1.17 of \cite{BK}.\end{proof}

Since $3$-Sasakian manifolds are abundant, the same is true for
generalized, Sasakian structures of bi-Sasakian type. For instance,
the unit spheres $S^{4n+3}$ have two distinct $3$-Sasakian
structures (Example 13.2.6 in \cite{BK}).

It seems to be difficult to find examples of non-bi-Sasakian type,
but, we can formulate the required conditions in the general case (a
different formulation was given in \cite{VCRF}).
\begin{theorem}\label{thSgen} A generalized Sasakian structure
on the manifold $M$ is equivalent with a pair
$(F_\pm,Z_\pm,\xi_\pm,\gamma)$ of classical, normal, almost
contact, metric structures that satisfy the following conditions
\begin{equation}\label{semnederivate}
L_{Z_+}\Xi_+=-L_{Z_-}\Xi_-,\end{equation}
\begin{equation}\label{bGualt61}
\Xi_\pm-d\xi_\pm+L_{Z_\pm}L_{Z_\pm}\Xi_\pm=0,\end{equation}
\begin{equation}\label{bGualt71}  d\Xi_\pm-\xi_\pm
\wedge L_{Z_\pm}\Xi_\pm+(dL_{Z_\pm}\Xi_\pm)^c=0,\end{equation}
supplemented by a $1$-form $\kappa\in\Omega(M)$.\end{theorem}
\begin{proof}
The upper index $c$ in (\ref{bGualt71}) comes from the following
notation inspired by complex geometry:
$\forall\lambda\in\Omega^k(M)$, $\lambda^c$ is the form given by
$$\lambda^c(X_1,...,X_k)=\lambda(F_\pm X_1,...,F_\pm
X_k).$$

Let $(\tilde{G},\mathcal{J})$ be the corresponding, generalized
K\"ahler structure on $M\times\mathds{R}$ with the corresponding
structures $(\Gamma,\Psi,J_\pm)$. The latter are equivalent with a
pair of normal, almost contact metric structures
$(F_\pm,Z_\pm,\xi_\pm,\gamma)$ on $M$.

Furthermore, the K\"ahler forms (\ref{formaomega}) of
$(\Gamma,J_\pm)$ have to satisfy the characteristic conditions
(\ref{eqluiGaltieri}) of a generalized K\"ahler structure. These
conditions are
\begin{equation}\label{bGualt} \begin{array}{l}d\omega_\pm(J_\pm
(X+a\frac{\partial}{\partial t}), J_\pm
(Y+b\frac{\partial}{\partial t}),J_\pm(U+u\frac{\partial}{\partial
t})\vspace{2mm}\\=\pm d\Psi(X+a\frac{\partial}{\partial
t},Y+b\frac{\partial}{\partial t},U+u\frac{\partial}{\partial
t}),\end{array}\end{equation} equivalently (use (\ref{Jacontact}),
(\ref{formaomega}), (\ref{bGualt})),
\begin{equation}\label{bGualt1}
\begin{array}{l} d\Xi_\pm(F_\pm X,F_\pm Y,F_\pm U)
+\sum_{Cycl}u[i(Z_\pm)d\Xi_\pm](F_\pm X,F_\pm Y)\vspace{2mm}\\
-\sum_{Cycl}\xi_\pm(U)(\Xi_\pm-d\xi_\pm)(F_\pm X+aZ_\pm,F_\pm Y
+bZ_\pm)\vspace{2mm}\\ =\pm\{d\psi(X,Y,U)+
\sum_{Cycl}u(\psi+d\kappa)(X,Y)\},
\end{array}\end{equation} where the cyclic permutations
are on the arguments $(X,a),(Y,b),(U,u)$.

Furthermore, (\ref{bGualt}) may be decomposed into the following
two cases 1) $a=0,b=0,u=1,U=0$ and 2) $a=b=u=0$. In case 1),
(\ref{bGualt1}) reduces to
\begin{equation}\label{bGualt2}\begin{array}{c} [i(Z_\pm)d\Xi_\pm](F_\pm X,F_\pm
Y)+\{\xi_\pm\wedge[i(Z_\pm)(\Xi_\pm-d\xi_\pm)]\circ F_\pm\}(X,Y)
\vspace{2mm}\\
=\pm(\psi+d\kappa)(X,Y),\end{array}\end{equation} which, by taking
into account $i(Z_\pm)\Xi_\pm=0$, $\xi_\pm(Z_\pm)=1$,
$i(Z_\pm)d\xi_\pm=L_{Z_\pm}\xi_\pm=0$ (because of normality),
becomes
\begin{equation}\label{bGualt20}[L_{Z_\pm}\Xi_\pm](F_\pm X,F_\pm
Y)=\pm(\psi+d\kappa)(X,Y).\end{equation} But, $\Xi_\pm(F_\pm
X,F_\pm Y)=\Xi_\pm(X,Y)$ and, by normality, $L_{Z_\pm}F_\pm=0$;
thus, we see that (\ref{bGualt20}) may be written under the form
\begin{equation}\label{valoarepsi}
\psi+d\kappa=L_{Z_+}\Xi_+=-L_{Z_-}\Xi_-.\end{equation} Formula
(\ref{valoarepsi}) yields (\ref{semnederivate}) and defines $\psi$
if $kappa$ is given.

In case 2), (\ref{bGualt1}) reduces to
\begin{equation}\label{bGualt3} \begin{array}{c}
d\Xi_\pm(F_\pm X,F_\pm Y,F_\pm U)
-\sum_{Cycl}\xi_\pm(U)[(\Xi_\pm-d\xi_\pm)(F_\pm X,F_\pm
Y)]\vspace{2mm}\\ =\pm d\psi(X,Y,U).\end{array}\end{equation}
Since $TM=im\,F_\pm\oplus span\{Z_\pm\}$, (\ref{bGualt3}) is
equivalent with the pair of conditions where $(X,Y,U)$ are taken:
(i) $(F_\pm X,F_\pm Y,Z_\pm)$, (ii) $(F_\pm X,F_\pm Y,F_\pm U)$.
In case (i) (\ref{bGualt3}) becomes
\begin{equation}\label{bGualt6}
\Xi_\pm-d\xi_\pm=\mp [i(Z_\pm)d\psi]^c,\end{equation} and
in case (ii) (\ref{bGualt3}) becomes
\begin{equation}\label{bGualt7}
d\Xi_\pm-\xi_\pm\wedge(i(Z_\pm)d\Xi_\pm) =\mp
(d\psi)^c.\end{equation}

Thus, the conditions that characterize the generalized Sasakian
case are (\ref{valoarepsi}), (\ref{bGualt6}) and (\ref{bGualt7}).
Then, since (\ref{valoarepsi}) gives
$d\psi=dL_{Z_+}\Xi_+=-dL_{Z_-}\Xi_-$, we may replace
(\ref{bGualt6}) and (\ref{bGualt7}), respectively, by
\begin{equation}\label{bGualt60}
\Xi_\pm-d\xi_\pm=-[i(Z_\pm)d(L_{Z_\pm}\Xi_\pm)]^c,\end{equation}
\begin{equation}\label{bGualt72}  d\Xi_\pm-\xi_\pm
\wedge(i(Z_\pm)d\Xi_\pm)=-(dL_{Z_\pm}\Xi_\pm)^c.\end{equation}
Above, we may replace $i(Z_\pm)d\Xi_\pm=L_{Z_\pm}\Xi_\pm$ and
$i(Z_\pm)d(L_{Z_\pm}\Xi_\pm)=L_{Z_\pm}L_{Z_\pm}\Xi$ (because of
$i(Z_\pm)\Xi_\pm=0$), and
$[L_{Z_\pm}L_{Z_\pm}\Xi]^c=L_{Z_\pm}L_{Z_\pm}\Xi$. The results
exactly are (\ref{bGualt61}) and (\ref{bGualt71}).
\end{proof}
\begin{rem}\label{obsremote} {\rm  The $2$-form
$L_Z\Xi$ might be called the {\it derived fundamental form}. On the
other hand, a classical, normal, almost contact, metric structure
$(F,Z,\xi,\gamma)$ such that
\begin{equation}\label{eqremoteS}\begin{array}{l}
\Xi-d\xi+L_ZL_Z\Xi=0,\vspace{2mm}\\
d\Xi-\xi \wedge L_{Z}\Xi+(dL_{Z}\Xi)^c=0\end{array}\end{equation}
might be called a {\it remotely Sasakian structure}. The last term
of the left hand side of the second condition (\ref{eqremoteS}) is
also equal to $L_Z(d\Xi)^c$. It follows easily that a Sasakian
structure is remotely Sasakian and a remotely Sasakian structure
that satisfies the condition $L_ZL_Z\Xi=0$ is Sasakian. With this
terminology, a generalized Sasakian structure is equivalent with a
pair $(F_\pm,Z_\pm,\xi_\pm,\gamma)$ of remotely Sasakian
structures, with the same metric, and with sign-opposite derived
fundamental forms, complemented by an arbitrary $1$-form $\kappa$.
Unfortunately, we do not have a real understanding of the
non-Sasakian, remotely Sasakian structures.}\end{rem}

We mention the following corollaries of Theorem \ref{thSgen}.
\begin{corol}\label{cazbiS} A pair $(F_\pm,Z_\pm,\xi_\pm,\gamma)$
of normal, almost contact, metric structures, with vanishing
derived fundamental forms (equivalently, with $L_{Z_\pm}\gamma=0$)
corresponds to a generalized, Sasakian structure iff the
structures are Sasakian.\end{corol}
\begin{proof} If there is a corresponding generalized, Sasakian
structure, the conclusion follows from (\ref{bGualt61}).
Conversely, if the structures are Sasakian, the left hand sides of
(\ref{bGualt61}), (\ref{bGualt71}) vanish and we get a
generalized, Sasakian structure by adding $\psi=-d\kappa$ for an
arbitrary $1$-form $\kappa$.\end{proof}
\begin{corol}\label{concluzie2} If the form
$\psi$ of a generalized Sasakian manifold is a closed $2$-form, the
corresponding structures $(F_\pm,Z_\pm,\xi_\pm,\gamma)$ $M$ are
classical Sasakian structures and $\psi$ must be an exact
form.\end{corol}
\begin{proof} If $d\psi=0$, then  (\ref{valoarepsi}) yields
$dL_{Z_\pm}\Xi_\pm=0$ and the conclusion follows from
(\ref{bGualt61}).\end{proof}
\begin{rem}\label{peculiar} {\rm There is a peculiarity in the
terminology that we have chosen: a generalized Sasakian structure
in the sense of Definition \ref{defgenS} may not be a generalized,
almost contact structure. Our terminology is motivated by the
equivalence with a pair of classical, normal, almost contact
metric structures $(F_\pm,Z_\pm,\xi_\pm,\gamma)$ on $M$
(satisfying some supplementary conditions). Such a pair provides
the pair $J_\pm$ of translation invariant, complex structures on
$M\times\mathds{R}$. But, if we use $(\Gamma,\Psi,J_\pm)$ to
reconstruct the generalized structure $
\mathcal{J}$ (using formulas (\ref{AdinJ})), we get a translation
invariant structure that may not be $M$-adapted.}\end{rem}
 %\end{center}
\hspace*{7.5cm}{\small \begin{tabular}{l} Department of
Mathematics\\ University of Haifa, Israel\\ E-mail:
vaisman@math.haifa.ac.il \end{tabular}}
\end{document}